%
%
%

\let\he=\heads
\let\hee=\headss

\def\hh #1\par{\underbar{\he#1}\vskip1pc}
\def\hhh #1\par{\underbar{\hee#1}\vskip1.2pc}
\def\ar #1\ {\bf#1}
\font\title=cmbx12
\font\by=cmr8
\font\author=cmr10
\font\adress=cmsl10
\font\abstract=cmr8

\def\tit #1\par{\centerline{\title #1}}
\def\bby{\centerline{\by BY}\par}
\def\aut #1\par{\centerline{\author #1}}
\def\abs{\abstract\centerline{ABSTRACT}\par}
\def\abss #1\par{\abstract\midinsert\narrower\narrower\noindent #1\endinsert}
\def\\{\noindent}
\def\vv{\par\rm}
\def\vvv{\medbreak\rm}

\def\sec{\vskip 2pc}

\def\sh{[Sh]}
\def\gu{[GuSh]}

\def\qed{\line{\hfill$\heartsuit$}\medbreak}
\def\pt #1. {\medbreak\noindent{\bf#1. }\enspace\sl}
\def\pd #1. {\medbreak \noindent{\bf#1. }\enspace}
\def\dd{Definition\ } 
\def\bb #1{{\bar #1}}
\def\bbp{\bar P}
\def\fo{\ ||\kern -4pt-}
\def\om{\omega}
\def\al{\alpha}
\def\be{\beta}
\def\ga{\gamma}
\def\Ga{\Gamma}

\def\de{\delta}
\def\et{\eta}
\def\si{\sigma}
\def\ta{\tau}
\def\ale{\aleph_0}
\def\sb{\subseteq}

\def\imp{\Rightarrow}
\def\lan{\langle}
\def\ran{\rangle}

\def\all{\forall}
\def\ex{\exists}
\def\qmu{$Q$\llap{\lower5pt\hbox{$\sim$}\kern1pt}$_\mu$}

\def\Z{\Bbb Z}
\def\Q{\Bbb Q}
\def\M{{\cal M}}
\def\em{\emptyset}
\def\too{\hookrightarrow}
\def\red{\restriction} 

\def\st{such that\ }

\def\int{interpretation\ }

\def\hd{Hdeg}
\def\today{\ifcase\month\or January\or February\or March\or April\or May\fi 
\space\number\day, \number\year}
\def\omn{\thinspace^{n\ge}\omega}
\def\ee #1,#2{{#1}\thinspace^\wedge\lan{#2}\ran}
\def\seg{\triangleleft}
\def\bina{\thinspace^{\omega>}2}

%
%
%
%
%
%
\baselineskip 1.2pc
\hsize=6in
\input mssymb
%
%
%
\line{\hfill\today}
\vskip.5in
\tit Uniformization, choice functions and well orders in the class of trees.
\par
\vskip.5in
\bby
\aut SHMUEL LIFSCHES and SAHARON SHELAH\footnote*{The second author would like 
to thank the U.S.--Israel Binational Science Foundation for 
partially supporting this research. Publ. 539}
\par
{\adress Institute of Mathematics, The Hebrew University of Jerusalem, 
Jerusalem, Israel} \par
\vskip.5in
%
%
%
\abs
\abss The monadic second-order theory of trees allows quantification over 
elements and over arbitrary subsets. We classify the class of trees with 
respect to the question: does a tree $T$ have a definable choice function 
(by a monadic formula with parameters)? A natural dichotomy arises where the 
trees that fall in the first class don't have a definable choice function 
and the trees in the second class have even a definable well ordering of 
their elements. This has a close connection to the uniformization problem.
\par
\vskip.6in
%
%
\hh 0. Introduction  \par  \rm
The {\sl uniformization problem} for a theory $T$ in a language $L$ can be 
formulated as follows: Suppose 
$T\vdash(\all\bb Y)(\exists\bb X)\phi(\bb X,\bb Y)$ where $\phi$ is an 
$L$-formula and $\bb X,\bb Y$ are tuples of variables. Is there another 
$L$-formula $\phi^*$ \st
$$T\vdash(\all\bb Y)(\all\bb X)[\phi^*(\bb X,\bb Y)\imp\phi(\bb X,\bb Y)]
\ \ \ {\rm and} \ \ \ T\vdash(\all\bb Y)(\exists!\bb X)\phi^*(\bb X,\bb Y) ?$$
\\Here $\exists!$ means ``there is a unique''.

The monadic second-order logic is the fragment of the full second-order 
logic that allows quantification over elements and over monadic (unary) 
predicates only. The monadic version of a first-order language $L$ can be 
described as the augmentation of\/ $L$ by a list of quantifiable set variables 
and by new atomic formulas $t\in X$ where $t$ is a first order term and $X$ 
is a set variable. The monadic theory of a structure $\M$ is the theory of 
$\M$ in the extended language where the set variables range over all subsets 
of  $|\M|$ and $\in$ is the membership relation.

Given a tree $T$ we may ask the following question: is there a sequence 
$\bb P$ of subsets of\/ $T$ and a formula $\varphi(x,X,\bb Z)$ in the monadic 
language of trees \st

\\$T\models\varphi(a,A,\bb P)\ \imp \ [A\not=\em\ \&\ a\in A]$  \ \ 
$T\models(\all X)(\exists y)[X\not=\em\imp\varphi(y,X,\bb P)]$  \ and 

\\$T\models\varphi(a,A,\bb P)\wedge\varphi(b,A,\bb P)\ \imp a=b$\ ?

\\If the answer is positive we will say that $T$ has a (monadically) definable 
choice function (with parameters) and that $\varphi$ defines a choice function 
from non-empty subsets of\/ $T$. Note that if we let $\phi(x,Y)$ be the 
formula that says ``if\/ $Y$ is not empty then $x\in Y$'' then a negative 
answer to the choice function problem for $T$ implies a negative answer to the
uniformization problem for the monadic theory of\/ $T$ (with $\phi$ being a 
counter-example).

dealing with the choice function problem we split the class of trees into 
two natural parts, wild trees and tame trees and prove the following: 

\pt Theorem. Let $T$ be a tree. If\/ $T$ is wild or $T$ embeds $\bina$ then 
there is no definable choice function on $T$ (by a monadic formula with 
parameters). If\/ $T$ is tame and does not embed $\bina$ then there is even 
a definable well ordering of the elements of\/ $T$ by a monadic formula (with 
parameters) $\varphi(x,y,\bb P)$.
\vvv
\\Looking at the definitions and proofs we observe that a tree is tame [wild] 
if and only if it's completion is tame [wild] and that the counter-examples 
for the choice function problem are either anti-chains or linearily ordered 
subsets of\/ $T$. 
Hence we can prove:
\pt Conclusion. Let $T$ be a tree and $T'$ be it's completion. Then the 
following are equivalent:

\\a) For some $n,l<\om$, for every anti-chain/branch $A$ of\/ $T$ there is a 
monadic formula $\varphi_A(x,X,\bb P_A)$ with quantifier depth $\le n$ and 
$\le l$ parameters from $T$, that defines a choice function from non empty 
subsets of\/ $A$. 

\\b) There is a monadic formula, with parameters, $\psi(x,y,\bb P)$ that 
defines a well ordering of the elements of\/ $T$.

\\c) There is a monadic formula, with parameters, $\psi'(x,y,\bb P')$ that 
defines a well ordering of the elements of\/ $T'$.
\vvv

\\The paper continues the work by Gurevich-Shelah (\gu) who answered 
negatively a question by Rabin ([Ra]), by showing that the answer for the 
choice function problem is negative in $\bina$.

\\The `positive' results on the existence of a definable well ordering 
($\S\S$3,5) are elementary and do not require knowledge of monadic logic. The 
negative results ($\S\S$2,3,4) are based on understanding of some composition 
theorems that hold for the monadic theory of trees. These facts are collected 
in $\S$1.

\\More details and Historical background can be found in [Gu] and \gu.
\sec

%
%
\hh 1. Composition Theorems \par
In this section we will define partial theories and establish the technical 
tools that will be applied later. We will formalize composition theorems that 
will enable to compute the partial theory of a tree from partial theories of 
it's parts. Using such theorems enables to prove that if for example a dense 
chain does not have definable choice function then a tree with a dense branch
does not have a definable choice function.
\pd \dd 1.1. 
$(T,\seg)$ is a {\sl tree} if\/ $\seg$ is a partial order on $T$ and
for every $\et\in T$, $\{\nu:\nu\seg\et\}$ is linearily ordered by 
$\seg$. 

\\Note, a chain $(C,<)$ and even a set without structure $I$   
is a tree. \vv
\pd \dd 1.2.  Let $T$ be a tree

\\1. $X\sb T$ is a {\sl convex subset} if\/ $\et,\nu\in X$ and 
$\et\seg \si\seg \nu\in T$ implies $\si\in X$. If\/ $T$ is a chain we use the 
term a {\sl convex segment} or just a {\sl segment}.

\\2. $(S,\seg)$ is a {\sl subtree} of\/ $(T,\seg)$ if\/ $S\sb T$ and 
$S$ is a convex subset of\/ $T$.

\\3. $B\sb T$ is a {\sl sub-branch} of\/ $T$ if\/ $B$ is convex and  
$\seg$--linearily ordered. 

\\4. $B\sb T$ is a {\sl branch} of\/ $T$ if\/ $B$ is a maximal sub-branch of\/ $T$.   

\\5. $A\sb T$ is an {\sl initial segment} of\/ $T$ if\/ $A$ is a sub-branch 
that is $\seg$--downward closed. $\et$ is {\sl above} an initial segment $A$  
if $\nu\in A\imp \nu\seg\et$.

\\6. For $\et\in T$, $T_{\ge\et}$ is the sub-tree 
$(\{\nu\in T: \et\seg\nu\},\seg)$.  \ 
$T_{>\et}$ is the sub-tree $(T_{\ge\et}\setminus\{\et\}\ ,\ \seg)$. \ 
For $A\sb T$ an initial segment, $T_{\ge A}$ and $T_{>A}$ are defined 
naturally.

\\7. For $\et\in T$ we deote by $suc(\et)$ or $suc_T(\et)$ the set of 
$\seg$--immediate successors of\/ $\et$ (which may be empty). 

\\8. For $\et,\nu\in T$ we denote the {\sl intersection of\/ $\et$ and $\nu$ in 
$T$} by $\et\wedge\nu$. This may be a member of\/ $T$ or an initial segment of 
$T$, in any case the meaning of\/ $\et\wedge\nu\seg\si$ is natural and 
$\si\seg\et\wedge\nu$ [$\si\in\et\wedge\nu$] is used only when $\et\wedge\nu$ 
is an element [an initial segment].

\\9. If there is an $\et\in T$ that satisfies $(\all\nu\in T)[\et\seg\nu]$ we 
say that $T$ has a root and denote $\et$ by $root(T)$.  

\\10. $\et,\nu\in T$ are {\sl incomparable in} $T$ if neither $\et\seg\nu$ 
nor $\nu\seg\et$. $X\sb T$ is an {\sl anti-chain of} $T$ if\/ $X$ consists of 
pairwise incomparable elements of\/ $T$.

\\11. A {\sl gap} in $T$ is a pair $(A,B)$ where $A\cap B=\em$, $A\cup B$ is a
sub-branch, $A$ is an initial segment, 
(so $\et\in A,\nu\in B\imp \et\seg\nu$), $A$ without a $\seg$-maximal element,
$B$ without a $\seg$-minimal element, and for some $\si\in T$ for every 
$\et\in A$ and $\nu\in B$ we have $\et\seg\si\ \&\ \nu,\si$ are incomparable. 

\\12. {\sl Filling a gap} $(A,B)$ in $T$ is adding a node $\ta$ to $T$ \st 
$\et\in A\imp\et\seg\ta$, $\nu\in B\imp\ta\seg\nu$ and for every $\si$ as in  
(11) we have $\ta\seg\si$.
\vv
\pd \dd 1.3. The {\sl full binary tree} is the tree $(\bina\ ,\ \seg)$ where 
for sequences $\et,\nu\in\bina$, $\et\seg\nu$ means $\et$ is an initial 
segment of\/ $\nu$.  \vv
\pd \dd 1.4. The {\sl monadic language of trees} $L$ is the monadic version 
of the language of partial orders $\{\seg\}$. Usually $\seg$ means ``smaller 
than or equal'' but when we restrict ourselves to chains 
(linearily ordered sets) we use $<$ and $\le$.  For simplicity, we add to $L$  
the predicate $sing(X)$ saying ``$X$ is a singleton'' so that we can quantify 
only over subsets. Note that everything that is defined in 1.2 is definable 
in $L$. \vvv
Next we define, following \sh, the partial theories of a tree $T$. These are 
finite approximations of the monadic theory of\/ $T$. $Th^n(T;\bb P)$ is 
essentially the monadic theory of\/ $(T,\bb P,\seg)$ restricted to sentences of 
quantifier depth $n$.
\pd \dd 1.5. For any tree $T$, $\bb A\in {\cal P}(T)^{lg(\bb A)}$, 
and a natural number $n$, define by induction 
$$t = Th^n(T;\bb A).$$ 
for $n=0$: 
$$t =\bigl\{\phi(\bb X):\ \phi(\bb X)\in L,\ \phi(\bb X)
{\rm \ quantifier\ free},\ T\models\phi(\bb A) \bigr\}.$$
for $n=m+1$:
$$t = \bigl\{Th^m(T;\bb A\thinspace^\wedge B):\ B\in {\cal P}(T)\}.$$
$T_{n,l}$ is the set of all formally possible $Th^n(T;\bbp)$ where $T$ is a 
tree and $l(\bbp)=l$.
\vv
\pt Fact 1.6. (A) For every formula $\psi(\bb X)\in L$ there is an $n$  
such that from $Th^n(T;\bb A)$ we can effectively decide whether 
$T \models \psi(\bb X)$.  

\\(B) If\/ $m\ge n$ then $Th^n(T;\bb A)$ can be effectively computed 
from $Th^m(T;\bb A)$.  

\\(C) Each $Th^n(T;\bb A)$ is hereditarily finite, and we can effectively 
compute the set $T_{n,l}$ 
of formally possible $Th^n(T,\bb A)$. \vvv
Next we recall the composition theorem for linear orders which states that 
the partial theory of a chain can be computed from the partial theories of 
it's convex parts. This allows us to sum partial theories formally.
\pd \dd 1.7. If\/ $C,D$ are chains then $C+D$ \ is any chain that can be split
into an initial segment isomorphic to $C$ and a final segment isomorphic to 
$D$.  

\\If\/ $\langle C_i:i<\al\rangle$ is a sequence of chains then 
$\sum_{i<\al}C_i$ \ is any chain $D$ that is the concatenation of segments  
$D_i$, such that each $D_i$ \ is isomorphic to $C_i$. \vv
\pt Theorem 1.8 (composition theorem for linear orders). 

\\(1) If \ $l(\bb A)=l(\bb B)=l(\bb A')=l(\bb B')=l$, and
$$Th^m(C,\bb A) = Th^m(C',\bb A') \ \ {\rm and}\ \ 
Th^m(D,\bb B) = Th^m(D',\bb B')$$
then 
$$Th^m(C+D,A_0\cup B_0,\ldots,A_{l-1}\cup B_{l-1}) = 
Th^m(C'+D',A'_0\cup B'_0,\ldots,A'_{l-1}\cup B'_{l-1}).$$ 
(2) If \ $Th^m(C_i,\bb {A_i}) = Th^m(D_i,\bb {B_i})$, 
$l(\bb A_i)=l(\bb B_i)=l$ for each $i<\al$, then
$$Th^m\Bigl( \sum_{i<\al}C_i,\ \cup_i A_{1,i},\ldots,\cup_i A_{l-1,i}\Bigr) = 
Th^m\Bigl( \sum_{i<\al}D_i,\ \cup_i B_{1,i},\ldots,\cup_i B_{l-1,i}\Bigr).$$ 
\vv
\pd Proof. By \sh\ Theorem 2.4 (where a more general theorem is proved), 
or directly by induction on $m$. 
\vv \qed
\pd Notation 1.9. 

\\(1) $t_1+t_2=t_3$ means: for some $m,l<\om$, $t_1,t_2,t_3\in T_{m,l}$ 
(remember definition 1.5) and  \par
\\if 
$$t_1=Th^m(C,A_0,\ldots,A_{l-1}) \ {\rm and} \ 
t_2=Th^m(D,B_0,\ldots,B_{l-1})$$ 
then  
$$t_3=Th^m(C+D,A_0\cup B_0,\ldots,A_{l-1}\cup B_{l-1}).$$
By the previous theorem, the choice of\/ $C$ and $D$ is immaterial. 

\\(2) $\sum_{i<\al}Th^m(C_i,\bb {A_i})$ \ \ is  \ 
$Th^m(\sum_{i<\al}C_i,\ \cup_{i<\al}A_{1,i},\ldots,\cup_{i<\al} A_{l-1,i})$. 

\\(3) If\/ $D$ is a subchain of\/ $C$ and $X_1,\ldots,X_{l-1}$ \ are subsets of 
  $C$ then $Th^m(D,X_0,\ldots,X_{l-1})$ abbreviates 
  $Th^m(D,X_0\cap D,\ldots,X_{l-1}\cap D)$.  

\\(4) We use abbreviations as $\bb P\cup\bb Q$, $\cup_i\bb P_i$ and 
$\bb P\sb C$. The meanings should be clear. 

\\(5) For $C$ a chain, $a<b\in C$ and $\bb P\sb C$ we denote by 
$Th^n(C;\bb P)\red_{[a,b)}$ the theory $Th^n([a,b);\bb P\cap[a,b))$.

\vvv 
%
%
%
%
%
%
The class of trees has some weaker (but sufficient for our purpose) 
composition theorems. First we define the composition of subtrees of the full 
binary tree following \gu\ and quote the the respective composition theorem.
\pd \dd 1.10. Let $M\sb\bina$ be a tree. A {\sl grafting function} on $M$ is 
a function $g$ satisfying the following conditions:

\\(a) $Dom(g)\sb M\times\{0,1\}$,

\\(b) if\/ $(x,0)\in Dom(g)$ then $\ee x,0\not\in M$ and
if\/ $(x,1)\in Dom(g)$ then $\ee x,1\not\in M$,

\\(c) every value $g(x,d)$ of\/ $g$ ($d\in\{0,1\}$) is a tree $\sb\bina$.

\\A {\sl composition} of a tree $M$ and a grafting function $g$ is the tree 
$$M\cup\big\{\ee x,d\thinspace^\wedge y :(x,d)\in Dom(g), y\in g(x,d) \big\}$$. \vv
\pt Theorem 1.11 (composition theorem for binary trees). Let $M\sb\bina$ 
be a tree, $N\sb\bina$ be the composition of\/ $M$ and a grafting function $g$, 
$\bb X\sb N$ and $n<\om$. \underbar{Then}, there is $m=m(n)<\om$ (effectively 
computable from $n$) \st from 
$Th^m(M;\bb X,\bb L^g(n,\bb X),\bb R^g(n,\bb X))$ we 
can effectively compute $Th^m(N;\bb X)$ where \par
\\$L^g_t(n,\bb X):=\big\{x\in M:(x,0)\in Dom(g), Th^n(g(x,0),\bb X)=t \big\}$

\\$\bb L^g(n,\bb X):=\big\{L^g_t(n,\bb X): 
t \ {\rm a \ formally \ possible \ n-theory} \big\}$

\\and $\bb R^g(n,\bb X)$ is defined similary by replacing $L,0$ with $R,1$.\vv
\pd Proof. This is theorem 2 in $\S$2.3. of \gu. The language that is used 
there is different from our $L$ but all the mentioned symbols are 
monadically inter-definable (with some additional parameters) with our $\seg$,
(For example the relation ``$X$ is an immediate left successor of\/ $Y$'' is 
easily definable from $\seg$ and the parameter 
$A:=\{\et\in\bina: (\ex\nu\in\bina)[\et=\ee\nu,0]\}$). Thus the translation 
of \gu's proof is clear.
\vv \qed
The next three theorems allow us to compute a partial theory  $Th^n(T;\bb X)$ 
from  partial theories of sub-structures of\/ $T$. The proofs are by 
induction on $n$ noting that $Th^0(T;\bb P)$ can express only statements as 
$P_i\seg P_j$, $P_i\in P_j$ and $P_i=P_j$ and that $Th^{n+1}$ is a collection 
of\/ $n$-theories. Everything is basically the same as in the previous case and  
we will not elaborate beyond that.

\pt Theorem 1.12 (composition theorem for general successors). Let $T$ be 
a tree, $\bb X\sb T$ and $A\sb T$ an initial segment (i.e. linearily ordered 
by $\seg$ and downward closed). 

\\For every $x$ above $A$ ($x\not\in A$ and 
$y\in A \imp y\seg x$) denote by $T_{A,x}$ the sub-tree 

\\$\{y\in T: (\exists z)[z\seg x\ \&\ z\seg y\ \&\ z\ {\rm above\ } A]\}$.

\\We say that $x$ and $y$ are equivalent above $A$ if\/ $x$ and $y$ are above 
$A$ and $T_{A,x}=T_{A,y}$ (compare with definition 4.1), finally let 
$\{T_i: i\in I_A\}$ list the equivalence classes above $A$ 
(it's a disjoint union of sub-trees).

\\\underbar{Then} for every $n<\om$, there is $m=m(n)<\om$ (effectively 
computable from $n$) \st from 
$Th^m(T_{\le A};\bb X)$ and 
$Th^m(I_A;\bb P^A(n,\bb X))$ we can effectively 
compute $Th^n(T;\bb X)$ where \par
\\$T_{\le A}:=\{y\in T: y \ {\rm not\ above\ } A\}$ 

\\$P^A_t(n,\bb X):=\{i\in I_A: Th^n(T_i;\bb X)=t\}$, 

\\$\bb P^A(\bb X):=\{P^A_t(n,\bb X): 
t \ {\rm a \ formally \ possible \ n-theory} \}$,  

\\and $Th^m(I_A;\bb P^A(n,\bb X))$ is the $m$-theory of a set without 
sructure -- i.e. in the monadic language of equality.

\\(A natural case is when for some $y\in T$ we have $A=\{z:z\seg y\}$, \  
$\{x_i:i\in I\}=suc(y)$ and $T_i=T_{\ge x_i}$). \vv \qed
\pt Theorem 1.13 (composition theorem for branches). Let $T$ be a tree, 
$B\sb T$ a branch, $\bb X\sb T$ and $n<\om$. $(B',\seg)$ is the chain that is 
obtained by adding nodes to fill the gaps in $B$ -- remember 1.2(12), 
(so $B'$ is contained in the completion of\/ $B$). And let $T'$ be the tree 
obtained by replacing the branch $B$ by $B'$ 

\\\underbar{Then} there is $m=m(n)<\om$ (effectively computable from $n$) \st 
from $Th^m(B';\bb P^{B'}(n,\bb X))$ we can effectively compute 
$Th^n(T;\bb X)$ where 

\\for $\et\in B'$, \ $T'^{B'}_{\ge\et}:=T'_{\ge\et}\setminus B'$ 

\\$P^{B'}_t(n,\bb X):=\big\{\et\in B:Th^m(T'^{B'}_{\ge\et};\bb X)=t\big\}$

\\$P^{B'}_t(n,\bb X):=
\big\{P^{B'}_t(n,\bb X): t \ {\rm a \ formally \ possible \ n-theory} \big\}$.

\\Moreover, if\/ $\bb Y\sb B$ then from $Th^m(B';\bb P^{B'}(n,\bb X),\bb Y)$ 
we can effectively compute $Th^n(T;\bb X,\bb Y)$.
\vv  \qed
\pd Notations 1.14. For stating the next composition theorem we need a 
considerable amount of notations.

\\Let $T$ be a tree, by ``$F\colon\bina\too T$ is an embedding'' we mean 
$F$ is 1-1 and for $\et,\nu\in\bina$, \ $\et\seg\nu\iff F(\et)\seg F(\nu)$,
we also assume that $T$ has a root and $F(root(\bina))=root(T)$.

\\let $S\sb T$ be $F''(\bina)$, it is a tree (but not necesarily a subtree of 
\/ $T$) that can be identified with $\bina$.

\\For $x=F(\et)\in S$ define $x^0\ [x^1] \in S$ to be $F(\ee\et,0)\ 
[F(\ee\et,1)]$. 

\\For $Y\sb S$ an anti-chain (hence an anti-chain of\/ $T$) let 
$Bush(Y):=\{x\in T:(\ex y\in Y)[x\seg y]\}$ (it's a subtree of\/ $T$) and let 
$Bush_S(Y):=Bush(Y)\cap S$ (it's a subtree of\/ $S$).

\\For every $y\in S$ denote $y^0\wedge y^1$ by $y^i$. It may be an element 
of\/ $T$ or an initial segment but remember the convention in 1.2(8).

\\For every $y\in S$ we define some subtrees of\/ $T_{\ge y}$ 
(some of them may be trivial if for example $y=y^i$):

\\0) $T_0(y)$:=$T_{\ge y}$.

\\1) $T_1(y)$:=$\big\{x\in T: (\neg y^i\seg x)\ \&\ 
(\ex z\not=y)[(z\seg x)\&(y\seg z\seg y^i)]\big\}$, [These are the elements    
that split from the segment $(y,y^i)$ ].

\\2) $T_2(y)$:=$\big\{x\in T: (y\seg x)\ \&\ 
(\all z)[(z\seg y^i)\&(z\seg x)\imp (z\seg y)]\big\}$. (If\/ $y^i$ is an 
initial segment replace $z\seg y^i$ with $z\in y^i$), [These are the elements 
that split from $y$ but not from the segment $(y,y^i)$ ].

\\3) $T_3(y)$:=$\big\{x\in T: (\neg y^0\seg x)\ \&\ 
(\ex z\not=y^i)[(z\seg x)\&(y^i\seg z\seg y^0)]\big\}$. 
(If\/ $y^i$ is an initial segment replace $(\ex z\not=y^i)$ with 
$(\ex z)(y^i\seg z)$ ), [These are the elements that split from the segment 
$(y^i,y^0)$ ].

\\4) $T_4(y)$:=$\big\{x\in T: (\neg y^1\seg x)\ \&\ 
(\ex z\not=y^i)[(z\seg x)\&(y^i\seg z\seg y^1)]\big\}$. 
(If\/ $y^i$ is an initial segment replace $(\ex z\not=y^i)$ with 
$(\ex z)(y^i\seg z)$ ), [These are the elements that split from the segment
$(y^i,y^1)$ ].

\\5) $T_5(y)$:=$\big\{x\in T: (y^i\seg x)\ \&\ 
(\all z)[(z\seg x)\&(z\seg y^0\vee z\seg y^1)\imp (z\seg y^i)]\big\}$. 
(If\/ $y^i$ is an initial segment replace $z\seg y^i$ with $z\in y^i$), 
[These are the elements that split from $y^i$ but not from the segments 
$(y^i,y^0)$ and $(y^i,y^1)$].

\\6) $T_6(y)$:= $T_{\ge{y^0}}$.

\\7) $T_7(y)$:= $T_{\ge{y^1}}$.

\\For $y\in S$, $\bbp\sb T$, $\bb t=\lan t_0,t_1\ldots,t_7\ran$, 
$t_i$ a possible $n$-theory, we have 
$y\in Q_{\bb t}\iff Th^n(T_0(y);\bbp)=t_0\ \&\ldots\&\ Th^n(T_7(y);\bbp)=t_7$.
For $y\not\in S$ we have $y\in Q_\em$.

\\Finaly let $\bb Q(n,\bbp)$ be $\lan Q_{\bb t}: \bb t$\ a possible sequence 
of\/ $n$-theories$\ran\thinspace^\wedge\lan Q_\em\ran$

\\Note that every anti-chain $Y$ is definable from $Bush_S(Y)$ and $S$ is 
definable from $\bb Q$.   \vv

\pt Theorem 1.15 (composition theorem for embeddings). Following the above 
notations, let $T$ be a tree and $F\colon\bina\too T$ an embedding.  

\\\underbar{Then} for every $Y\sb S$ an anti-chain, $y\in Y$, $\bb P\sb T$ and 
$n<\om$, there is $m=m(n)<\om$ (effectively computable from $n$) \st from 
$Th^m(Bush_S(Y);y,\bb Q(n,\bbp))$ 
we can effectively compute $Th^n(T;y,Y,\bbp)$.  
\vv \qed
\sec
%
%
%
%
%
\hh 2. Dense linear orders \par
Every finite set $A$ has a definable well ordering by a formula with $|A|$
parameters. This is not the case for infinite models. \vv
\pt Claim 2.1. Let $A$ be an infinite set without structure. Then there is no 
definable choice function on $A$. Moreover, if\/ $|A|>2^l$ then no formula with 
$\le l$ parameters defines a choice function on $A$.
\vv
\pd Proof. Let $\bb P\sb A$ and suppose $\varphi(x,X,\bb P)$ defines a choice 
function on an infinite $A$. Let $B\sb A$ be an indiscernible set with 
respect to (belonging to) $\bb P$ of size $\ge 2$. Then, for every 
$b_1, b_2 \in B$, $A\models \varphi(b_1,B,\bb P)$ 
iff\/ $A\models \varphi(b_2,B,\bb P)$, a contradiction.  
The second part is clear. \vv \qed
%
%
%
%
A chain $C$ that embeds a dense linear order (hence the rational order $\Q$) 
does not have a definable choice function. The proof is by applying a 
Ramsey-like theorem for additive colourings from \sh.
\pd \dd 2.2. 
\\(a) A {\sl colouring} of a chain $C$ is a function $f$ from the set of 
unordered pairs of distinct elements of\/ $C$, into a finite set $I$ of colours.

\\(b) The colouring $f$ is {\sl additive} if for $x_i<y_i<z_i\in C$ ($i=1,2$),
$$[f(x_1,y_1)=f(x_2,y_2),\ f(y_1,z_1)=f(y_2,z_2)] \ \imp \ 
f(x_1,z_1)=f(x_2,z_2).$$
\\In this case a partial operation $+$ is defined on $I$, \st for 
$x<y<z\in C$, $f(x,z)=f(x,y)=f(y,z)$. (Compare with 1.9(1)).

\\(c) A subchain $D\sb C$  is {\sl homogeneous} (for $f$) if there is an 
$i_0\in I$ \st for every $x<y\in D$, $f(x,y)=i_0$. \vv
\pt Theorem 2.3. If\/ $f$ is an additive colouring of a dense chain $C$, by a 
finite set $I$ of colours, then there is an interval of\/ $C$ which has a dense 
homogeneous subset. \vv
\pd Proof. This is theorem 1.3. in \sh.  \vv \qed
\pt Claim 2.4. Let $(C,<)$ be a linear order that embeds a dense linear order. 
Then there is no definable choice function on $C$. \vv
\pd Proof.  Let $\bb P\sb C$ and suppose $\varphi(x,X,\bb P)$ defines a choice 
function on $C$. Let $n$ be so that from $Th^n(C;x,X,\bb P)$ we know if 
$\varphi(x,X,\bb P)$ holds and finaly let $D\sb C$ be dense (in itself).
By 2.3 there is an $A\sb D$, dense inside an interval of\/ $D$, hence in itself, 
homogeneous with respect to the colouring $f(a,b)=Th^n(C;\bb P)\red_{[a,b)}$,
(Remember the notation 1.9(5)).

\\Let $t^*$ be the constant theory $Th^n(C;\bb P)\red_{[a,b)}$ for every $a<b$ 
in $A$. 
Let $\Z$ be the set of integers and $X\sb A$, $X:= \{x_n:n\in \Z\}$ be of 
order type $\Z$. Suppose our choice function picks $x_m$ from $X$, i.e. 
$C\models \varphi(x_m,X,\bb P)$. \par
\\We assume for simplicity of notations that $inf(X)$ and $sup(X)$ belong to 
$C$ and denote $inf(X)$ by $0$ and $sup(X)$ by $1$. So 
$Th^n(C;\bb P)\red_X = Th^n(C;\bb P)\red_{(0,1)}$. 

\\Letting $t_0$ be $Th^n(C;\bb P)\red_{\{x:x\le 0\}}$, and $t_1$ be 
$Th^n(C;\bb P)\red_{\{x:x\ge 1\}}$ we get:
$$ Th^n(C;\bb P) = t_0 + \sum_{k\in\Z}Th^n(C;\bb P)\red_{[x_k,x_{k+1})} + t_1
  =  t_0 + \sum_{k\in\Z}t^* + t_1$$
\\Now denote: \par
\\$t'_0$ := $Th^n(C;x_m,X,\bb P)\red_{\{x:x\le 0\}}$ 
( = $Th^n(C;\em,\em,\bb P)\red_{\{x:x\le 0\}}$ ), \par
\\$t'_1$ := $Th^n(C;x_m,X,\bb P)\red_{\{x:x\ge 1\}}$ 
( = $Th^n(C;\em,\em,\bb P)\red_{\{x:x\ge 1\}}$ ), \par
\\$t'$ := $Th^n(C;x_l,X,\bb P)\red_{[x_k,x_{k+1})}$  for $k\not=l$, 
( = $Th^n(C;\em,x_k,\bb P)\red_{[x_k,x_{k+1})}$ ) and \par 
\\$t^{(l)}$ := $Th^n(C;x_l,X,\bb P)\red_{[x_l,x_{l+1})}$ 
( = $Th^n(C;x_l,x_l,\bb P)\red_{[x_l,x_{l+1})}$ ).  \par 
\\Clearly $t_0$ determines $t'_0$, $t_1$ determines $t'_1$,  
$t'_0$ and $t'_1$ do not depend on $m$ and 
$t^*$ determines $t'$ and $t^{(l)}$. We also have, for every $l\in\Z$:
$$ Th^n(C;x_l,X,\bb P) = t'_0 + \sum_{j\in\Z,j<l}t' + t^{(l)} + 
\sum_{j\in\Z,j>l}t' + t_1$$
\\But, by homogeneity, we get for every $k,l\in\Z$: 

\\1) $t^{(k)}=t^{(l)}$, 

\\2) $Th^n(C;x_l,X,\bb P)\red_{(0,x_l)} =  \sum_{j\in\Z,j<l}t' = 
\sum_{j\in\Z,j<k}t' = Th^n(C;x_k,X,\bb P)\red_{(0,x_k)}$, 

\\3) $Th^n(C;x_l,X,\bb P)\red_{(x_l,1)} =  \sum_{j\in\Z,j>l}t' = 
\sum_{j\in\Z,j>k}t' = Th^n(C;x_k,X,\bb P)\red_{(x_k,1)}$. 

\\It follows that $Th^n(C;x_m,X,\bb P) = Th^n(C;x_l,X,\bb P)$ for every 
$l\in\Z$, but $\varphi$ ``chooses'' $x_m$ from $X$, (and can be computed from 
$Th^n$) -- a contradiction. \vv \qed
\sec
%
%
\hh 3. Scattered orders \par
A scattered order is a linear order that does not embed a dense order. 
We will define \hd, the Hausdorff degree of scattered chains, and show 
that a scattered chain $(C,<_C)$ has a definable well ordering if 
$\hd(C)<\om$ and that $\hd(C)\ge\om \ \imp $ there is no definable choice 
function on $C$. \vv
\pd \dd 3.1. We define by recursion the Hausdorff degree of a scattered chain 
$(C,<_C)$:

\\$\hd(C)=0$ iff\/ $C$ is finite \par
\\$\hd(C)=\al$ iff\/ $\wedge_{\be<\al}\hd(C)\not=\be$ and 
$C=\sum_{i\in I}C_i$ where $I$ is well ordered or inversely well ordered and 
for every $i\in I$, $\vee_{\be<\al}\hd(C_i)=\be$. \par
\\$\hd(C)\ge\de$ iff\/ $(\all\al<\de)(\hd(C)>\al)$ ($\de$ limit). \vv
\pt Claim 3.2. (1) Let $C$ be a scattered chain with $\hd(C)=\al$, 
$C'$ the completion of\/ $C$ and $D\sb C'$. Then $C'$ and $D$ are scattered and  
$\hd(D)\le\hd(C')=\al$.

\\(2) Let $C$ be a scattered chain. $\hd(C)$ is well defined 
(i.e. it is an ordinal $\al$).
\vv
\pd Proof. (1) By induction on $\al$.

\\(2) By [Ha].
\vv \qed
\pt Claim 3.3. Let $C$ be a scattered chain with $\hd(C)=n$. Then there are 
$\bb P\sb C$, $lg(\bb P) = n-1$, and a formula (depending on $n$ only)  
$\varphi_n(x,y,\bb P)$ that defines a well ordering of\/ $C$.  \vv
\pd Proof. By induction on $n=\hd(C)$: 

\\$n\le 1$: \ $\hd(C)\le 1$ implies $(C,<_C)$ is well ordered or inversely 
well ordered. A well ordering of\/ $C$ is easily definable from $<_C$.

\\$\hd(C)=n+1$: Suppose $C=\sum_{i\in I}C_i$ and each 
$C_i$ is of Hausdorff degree $n$. By the induction hypothesis there are a 
formula $\varphi_n(x,y,\bb Z)$ and a sequence $\lan \bb P^i: i\in I \ran$\  
with $\bb P^i\sb C_i$, $\bb P^i = \lan P_1^i,\ldots,P_{n-1}^i \ran$ \st
$\varphi_n(x,y,\bb P^i)$ defines a well ordering of\/ $C_i$. 

\\Let for $0<k<n$, $P_k := \cup_{i\in I}P^i_k$  (we may assume that the union 
is disjoint) and $P_n := \cup\{C_i : i {\rm \ even}\}$.

\\We will define an equivalence relation $\sim$ by $x\sim y$ iff 
$\bigwedge_i (x\in C_i \Leftrightarrow y\in C_i)$. 

\\$\sim$ and $[x]$, (the equivalence class of an element $x$), are easily 
definable from $P_n$ and $<_C$. We can also decide from $P_n$ if\/ $I$ is well 
or inversely well ordered (by looking at subsets of\/ $C$ consisted of 
nonequivalent elements) and define $<'$ to be $<$ if\/ $I$ is well ordered 
and the inverse of\/ $<$ if not.
$\varphi_{n+1}(x,y,P_1,\ldots,P_n)$ will be 
defined by: 
$$\varphi_{n+1}(x,y,\bb P) \Leftrightarrow 
\big[ x\not\sim y \ \&\  x<'y\big] \vee 
\big[x\sim y \ \&\ \varphi_{n}(x,y,P_1\cap[x],\ldots,P_{n-1}\cap[x])\big]$$ 

\\$\varphi_{n+1}(x,y,\bb P)$ well orders $C$. \vv \qed
Next we prove that a scattered orders of infinite $\hd$ don't have a 
definable choice function (hence a well ordering). 

\pd \dd 3.4. We define for every $n<\om$ a model $\M^n$ in the language 
consisted of a binary relation $<^n$: 

\\a) The universe of\/ $\M^n$, which will be denoted by $M^n$, is the tree 
$\omn$.   

\\b) Let, for every $\eta\in \omn$, $<_\eta$ be a linear ordering of 
$suc(\eta) := \{\eta^\wedge\lan k\ran : k<\om\}$ \st if\/ $lev(\eta)$ is even 
then $k<l\imp \eta^\wedge\lan k\ran <_\eta \eta^\wedge\lan l\ran$, 
and if\/ $lev(\eta)$ is odd then 
$k<l\imp \eta^\wedge\lan l\ran <_\eta \eta^\wedge\lan k\ran$. \par
\\(So $<_\eta$ orders $suc(\eta)$ with order type $\om$ if\/ $\eta$ is in an 
even level and with order type $\om^*$ if\/ $\eta$ is in an odd level). 

\\c) $<^n$ is the lexicographic order induced by the orders $<_\et$ of 
immediate successors.

\\$(M^n,<^n)$ is hence a chain. Note, the `usual' partial order $\seg$ on 
$\omn$ (being an initial segment), is not definable in $\M^n$. \vv
\pd \dd 3.5. We define by induction the scattered chains $C_n$ and $C_n^*$: 

$C_1 := \om$, \ \ $C_1^* := \om^*$,    \par
$C_2 := \sum_{i\in\om}\om^*$, \ \ $C_2^* := \sum_{i\in\om^*}\om$,  \par
\\and in general:  \par
$C_n := \sum_{i\in\om}C_n^*$, \ \ $C_n^* := \sum_{i\in\om^*}C_n$.  \vv
\pd \dd 3.6. $f\colon\M^n\too C$ is an embedding of\/ $\M^n$ in a 
scatterd chain $(C,<_C)$ if 
$f$ is 1--1 and $\sigma<^n\tau \imp f(\sigma)<_C f(\tau)$  \par
\pt Fact 3.7. Let $C$ be a scattered chain with \ $\hd(C)\ge n+1$. Then 
there is an embedding $f\colon\M^n\too C$.
\pd Proof. Clearly the following hold:  

$(\al)$ \ For a scattered chain $C$: 
$\hd(C)=n\imp [C_n\sb C$ or $C^*_n\sb C]$.

$(\be)$ \ $\M^n\sb\M^{n+1}$ 

$(\ga)$ \ There is an embedding $g\colon\M^n\too C_n$.

\\Now assume $\hd(C)=n+1$ and use $(\al)$. In the case $C_{n+1}\sb C$ we 
have by $(\ga)$ an embedding $g\colon\M^{n+1}\too C$ and by $(\be)$ an 
embedding $f\colon\M^n\too C$. In the case $C^*_{n+1}\sb C$ we have, by 
the definition of\/ $C^*_{n+1}$, $C_n\sb C^*_{n+1}$ and by $(\ga)$ an 
embedding $f\colon\M^n\too C$. 
\vv \qed
\pt Conclusion 3.8. Let $C$ be a scattered chain with \ $\hd(C)\ge\om$. Then, 
for every $n<\om$ there is an embedding of\/ $\M^n$ into $C$. \vv \qed
\pt Lemma 3.9. If\/ $C$ is scattered and $\hd(C)\ge\om$ then no monadic formula 
\ $\varphi(x,X,\bb P)$ defines a choice function on $C$. \vv
\pd Proof. Assume towards a contradiction that there is $\bb P\sb C$,  
$lg(\bb P)=l$ and $\varphi(x,X,\bb P)$ defines a choice function on $C$.
Let $m$ be so that from $Th^m(C;x,X,\bb P)$ we can decide if 
$C\models\varphi(x,X,\bb P)$.  As in the proof of 2.4 it is enough
to find an $B\sb C$, of order type $\Z$, homogeneous with respect to the 
colouring $f(a,b)=Th^m(C;\bb P)\red_{[a,b)}$.  Let
$$ n > |\{Th^m(D;\bb Q): D {\rm \ a\ chain\ },\ \bb Q\sb D,\ l(\bb Q)=l\}| = 
|T_{m,l}| $$  
and $f\colon\M^n\too C$ be an embedding. Let $T\sb C$ be the image of\/ $f$ 
and we will identify $T$ with $\omn$ and the submodel 
$(T,<_C)\sb(C,<_C)$ with the model $(\omn,<^n)$ 
\medbreak
\\Notation: We will write $<$ instead of\/ $<_C$ and it's restriction $<^n$. 
Given $\et<\nu\in\omn=T$ we will write $Th^m[\eta,\nu)$ instead of
$Th^m(C;\bbp)\red_{[\eta,\nu)}$.  $T_{\ge\et}$ and $T_{>\et}$ are the usual 
subsets of\/ $\omn=T$
\medbreak 
\\We will begin to thin out the tree $T=\omn$, in order to obtain a quite 
homogeneous subtree $A\sb T$ going down with the levels. Arriving to a node 
$\et$, we will have defined $A_{\ge\nu}$ for every $\nu\in suc(\et)$ and 
will define  $A_{\ge\et}$ by thinning out $suc(\et)$ to a set $B_\et$ and 
taking $\{\et\}\cup\{A_{\le\nu}:\nu\in B_\et\}$. $A_{\ge\et}$ will satisfy 
the following:
$$ (*) \ \ [\si<\ta\in A_{\ge\et}, \ lev(\si)=lev(\ta)] \ \imp \ 
Th^m[\si,\ta) {\rm \ depends\ only\ on\ } lev(\si\wedge\ta)$$
Assume w.l.o.g that $n$ is odd.

\\\underbar{Step 1}: for every $\eta\in \omn$ with $lev(\eta)={n-1}$ pick 
out an infinite set $B_\eta\sb\om$ \st 
$$ k<l\in B_\eta \ \imp \ Th^m(C;\bbp)\red_{[\ee\et,k,\ee\et,l)}=t_\et $$ 
\\(note that $k<l<\om \imp \ee\et,k<\ee\et,l$), let $k_\et$ be the second 
element of\/ $B_\et$. 
Let $A_{\ge\et}$ be $\{\et\}\cup\{\ee\et,k : k_\et\le k\in B_\et\}$, this is 
a subtree of\/ $T$ and $(*)$ clearly holds.

\\\underbar{Step 2}: Given $\nu\in \omn$ with $lev(\nu)={n-2}$ we have 
defined $B_\si$, $k_\si$ and $A_{\ge\si}$ for every $\si\in suc(\nu)$. 
Pick out an infinite $B^0_\nu\sb\om$ so that $(*)$ will hold for $lev(\si)=
lev(\ta)=n-1$, $lev(\si\wedge\ta)=n-2$ \ i.e.
$$ k>l\in B^0_\nu \ \imp \ Th^m[\ee\nu,k,\ee\nu,l)=t_\nu $$ 
\\($suc(\nu)$ are ordered as $\om^*$).  Thin out $B^0_\nu$ to an infinite 
$B^1_\nu\sb\om$ so that $(*)$ will hold for $lev(\si)=lev(\ta)=n$, 
$lev(\si\wedge\ta)=n-2$ \ i.e.
$$ k>l\in B^1_\nu,\ \si=\ee\nu,k,\ \ta=\ee\nu,l\  \imp \ 
Th^m[\ee\si,{k_\si},\ee\ta,{k_\ta}) {\rm \ is \ constant} $$
\\Why does it suffice to look only at e.g. $\ee\si,{k_\si}$ ? because by the 
choice of\/ $t_\si$ and $A_\si$ we have $t_\si+t_\si=t_\si$ hence for every 
$\et<\si\in T$ we can break the paths $[\et,\ee\si,{k_\si})$ and 
$[\et,\ee\si,l)$, for $l\in A_\si$, into three parts: 
first from $\et$ to $\si$ then from $\si$ to it's `first' successor in 
$A_\si$, and then to $\ee\si,{k_\si}$ or $\ee\si,l$ 
(this is why we chose $k_\si$ to be the second element of\/ $A_\si$), 
but adding the last theory does not change the sum hence 
$Th^m[\et,\ee\si,{k_\si})=Th^m[\et,\ee\si,l)$ for every $l\in A_\si$. 
By a similar argument we can show that for every $l\in A_\si$ we have
$Th^m[\ee\si,{k_\si},\et)=Th^m[\ee\si,l,\et)$. 

\\Next, thin out $B_\nu^1$ to get $B_\nu$ so that $(*)$ will hold for 
$lev(\si)=lev(\ta)=n$, $lev(\si\wedge\ta)=n-1$ \ i.e.
$$k\in B_\nu,\ \si=\ee\nu,k \ \imp \ t_\si {\rm \ is \ constant} $$
\\let $k_\nu$ to be the second element of\/ $B_\nu$.  Define the subtree 
$A_{\ge\nu}$ to be $\{\nu\}\cup\{A_{\ee\nu,k} : k_\nu\le k\in B_\nu\}$.
Clearly $A_{\ge\nu}$ satisfies $(*)$.

\\\underbar{Step $n-1$}: we have reached $\lan e \ran$, the root of\/ $\omn$. 
\  $B^0_e, B^1_e, \ldots , B^{(n-1)\cdot(n-2)}_e=B_e$ are defined as before, 
taking care of $(*)$ for all the possibilities of the form 
$lev(\si)=lev(\ta)=k$, $lev(\si\wedge\ta)=l$ (some thinning outs are not 
necessary as they have been taken care of in previous steps), $k_e$, $t_e$ 
and $A_{\ge e}=A$ are defined as well.

\\\underbar{Final Step}:  By our construction, for every $\et<\nu$ in $A$, 
with $lev(\et)=lev(\nu)$, \ $Th^m(C;\bb P)\red_{[\et,\nu)}$ depends 
only on $lev(\et\wedge\nu)$  and we define $t_k$ by: 
$$t_k := Th^m[\et,\nu) {\rm\ where\ } \et<\nu, lev(\et)=lev(\nu)=n,\ 
lev(\et\wedge\nu)=n-k$$
By our choice of\/ $n$ we have some $k<l\le n$ with $t_k=t_l$. Let's show how 
to get a suitable homogeneous subset $B$ of\/ $T$ ($C$) from this.

\\Example 1. \ $t_1=t_2$  \par
\\Pick $\et\in A$ with $lev(\et)=n-2$. The successors of\/ $\et$ in $A$ have 
order type $\om^*$ and for every successor of\/ $\et$ in $A$, it's successors 
have order type $\om$.  Define:
$$B_1:= \big\{\ee\et,l\thinspace^\wedge\lan k_{\ee\et,l}\ran \ : \ 
l\in A_\et, l>k_\et \big\}$$
\\and
$$B_2:= \big\{\ee\et,{k_\et}\thinspace^\wedge\lan k\ran \ : \ 
k\in A_{\ee\et,{k_\et}} \big\}$$
\\and let $B=B_1\cup B_2$.

\\Clearly $B_1$ has order type $\om^*$, $B_2$ has order type $\om$ and $B$ 
has order type $\Z$. Moreover, for every $\si<\ta\in B_1$ we have 
$Th^m[\si,\ta)=t_1$ (since $lev(\si\wedge\ta)=n-1$) and for every 
$\si<\ta\in B$ with $\ta\in B_2$ we have $Th^m[\si,\ta)=t_2$ 
(since $lev(\si\wedge\ta)=n-2$).  By $t_1=t_2$ we conclude:
$$\all(\si<\ta\in B)\big[Th^m(C;\bb P)\red_{[\si,\ta)}=t_1\big].$$ 
\\Finding a homogeneous subset of\/ $C$ of order type $\Z$, we can proceed as 
in claim 2.2 to get a contradiction to ``$\varphi(x,X,\bb P)$ defines a choice 
function on $C$''.

\\Example 2. \ $t_2=t_3$  \par
\\Pick $\et\in A$ with $lev(\et)=n-3$. The successors of\/ $\et$ in $A$ have 
order type $\om$ and for every successor of\/ $\et$ in $A$, it's successors 
have order type $\om^*$.  Let $\si:=\ee\et,{k_\et}$ ($lev(\si)=n-2$), for 
$l>k_\et\in A_\et$ \ $\et_l:=\ee\et,l$ ($lev(\et_l)=n-2$) and 
$\si_l:=\ee{\et_l},{k_{\et_l}}$ ($lev(\si_l)=n-1$). Define:
$$B_1:= \big\{\ee{\si_l},{k_{\si_l}} \ : \ l\in A_\et, l>k_\et\big\}$$
\\($B_1$ has order type $\om^*$). 
\\To define $B_2$ we let, for $l\in A_\si$, $\ta_l:=\ee\si,l$ 
($\ta_l$ are extensions of\/ $\et$ and $\si$ with $lev(\ta_l)=n-1$) 
and then extend each $\ta_l$ to a $\rho_l$ defined by 
$\rho_l:=\ee{\ta_l},{k_{\ta_l}}$. So 
$$B_2:=\big\{\rho_l \ : \ l\in A_\si \big\}$$
\\and it has order type $\om$.  \ $B:=B_1\cup B_2$ has order type $\Z$ and we 
can easily check that for every $\nu_1<\nu_2\in B$ we have 
$Th^m[\nu_1,\nu_2)=t_3$ (as $\nu_1\wedge\nu_2=\et$ so 
$lev(\nu_1\wedge\nu_2)=n-3$) and for every $\nu_1<\nu_2\in B_2$ we have 
$Th^m[\nu_1,\nu_2)=t_2$ (as $\nu_1\wedge\nu_2=\si$ so 
$lev(\nu_1\wedge\nu_2)=n-2$).  By $t_2=t_3$ we conclude:
$$\all(\si<\ta\in B)\big[Th^m(C;\bb P)\red_{[\si,\ta)}=t_2\big].$$ 
\\and we proceed as before.
 
\\What we did in both examples can be described as follows: we fixed a node 
$\et\in A$ and a successor $\si$ of\/ $\et$, we extended the other successors of 
$\et$ and the successors of\/ $\si$ in a ``canonical'' way, ($\nu$ is extended 
to $\ee\nu,{k_\nu}$) to nodes of level $n$. The result is a homogeneous 
subset of\/ $C$ of order type $\Z$.

\\General case. $l+1<r$, $t_l=t_r$  \par
\\Let $\si,\ta\in A$ be \st $lev(\si)=lev(\ta)=n$ and $lev(\si\wedge\ta)=n-r$,
so $Th^m[\si,\ta)=t_r$. Then find $\rho\in A$ with $\si<\rho<\ta$, 
$lev(\rho)=n$, $lev(\si\wedge\rho)=n-(l+1)$ and $lev(\rho\wedge\ta)=n-r$. What   
we get is the following equation:
$$ t_r \ = \ Th^m[\si,\ta) \ = \ Th^m[\si,\rho)+Th^m[\rho,\ta) \ = \ 
t_{l+1}+t_r $$
\\but $t_r=t_l$ hence 
$$ (*) \ t_l=t_{l+1}+t_l $$
\\Imitate this computation: let $\si,\ta\in A$ be \st $lev(\si)=lev(\ta)=n$ 
and $lev(\si\wedge\ta)=n-(l+1)$,
so $Th^m[\si,\ta)=t_r$ and find $\rho\in A$ with $\si<\rho<\ta$, 
$lev(\rho)=n$, $lev(\si\wedge\rho)=n-(l+1)$ and $lev(\rho\wedge\ta)=n-l$. What   
we get is the following equation:
$$ t_{l+1} \ = \ Th^m[\si,\ta) \ = \ Th^m[\si,\rho)+Th^m[\rho,\ta) \ = \ 
t_{l+1}+t_l $$
\\hence 
$$ (**) \ t_{l+1}=t_{l+1}+t_l $$
\\Combining $(*)$ and $(**)$ we get $t_{l+1}=t_l$. Now proceed as in example 1
(if\/ $l$ is odd) or as in example 2 (if\/ $l$ is even) by taking ``canonical 
extensions'' of successors to get the required homogeneous subset $B$ of 
order type $\Z$. \vv \qed
\pt Conclusion 3.10. For every $m,l<\om$ there is an $n<\om$ \st if 
$C$ is a scattered chain and $\hd(C)\ge n+1$ then $C$ does not have a 
definable choice function by a formula with quantifier depth $\le m$ 
and with $\le l$ parameters. \vv
\pd Proof. Let $n$ be larger than \ 
$|T_{m,l}|$.
Now if\/ $\hd(C)\ge n+1$ then we can embed $\omn$ into $C$ and immitate the 
previous proof.
\vv \qed
\sec
%
%
%
%
%
\hh 4. Wild trees \par
Intuitively, wild trees are trees that have a large amount of splitting 
($4.2(1)(i)$) or have `complicated' branches ($4.2(1)(ii)(iii)$), the next 
two definitions state this formally. \ Wild trees don't have a definable 
choice function ($4.6$).

\pd \dd 4.1.  Let $(T,\seg)$ be a tree 

\\(1) If\/ $A$ is an initial segment of\/ $T$ then $top(A)$ is 
$\{x\in T: (\all t\in A)[t\seg x]\}$.  (It's a tree).

\\(2) Let $A$ be an initial segment of\/ $T$ then the binary relation 
$\sim^0_A$ on $T\setminus A$ is defined by 
$$x\sim^0_A y \ \iff \ (\all t\in A)[t\seg x\equiv t\seg y]$$
\\(It's an equivalence relation that says ``$x$ and $y$ `break' $A$ in the 
same place'').

\\(3) Let $A$ be an initial segment of\/ $T$ then the binary relation 
$\sim^1_A$ on $T\setminus A$ is defined by 
$$x\sim^1_A y \ \iff \ [x\sim^0_A y] \ \& \ 
(\exists z)[z\seg x \ \& \ z\seg y \ \& \ z\sim^0_A x]$$
\\(It's an equivalence relation that divides -- for every initial segment 
$B\sb A$ -- $top(B)/\sim^0_B$ into disjoint subtrees). \vv
\pd \dd 4.2. (1) A tree $T$ is called {\sl wild} if either \par 
$(i)$ $sup\big\{|top(A)/\sim^1_A| : A\sb T\ {\rm an\ initial\ segment}\big\}
\ge\ale$ \ \ or 

$(ii)$ There is a branch $B\sb T$ and an embedding $f\colon\Q\to B$ \ \ or

$(iii)$ All the branches of\/ $T$ are scattered linear orders but  
$sup\big\{\hd(B) : B\ {\rm a\ branch\ of\ } T\big\} \ge\om$.

\\(2) A tree $T$ is {\sl tame for} $(n^*,k^*)$ if 
the value in $(i)$ is $\le n^*$,  $(ii)$ does not hold and the value in 
$(iii)$ is $\le k^*$ 

\\(3) A tree $T$ is {\sl tame} if\/ $T$ is tame for $(n^*,k^*)$ for some 
$n^*,k^*\le\om$. \vv
\pt Claim 4.3. If\/ $T$ is a wild tree and (1)$(i)$ of 4.2 holds then no monadic 
formula \ $\varphi(x,X,\bb P)$ defines a choice function on $T$. \vv
\pd Proof. We will use the composition theorem for general successors 1.12. 

\\Suppose $\varphi(x,X,\bb P)$ defines a choice function on $T$ and 
$Th^n(T;x,X,\bb P)$ computes $\varphi$. For an initial segment $A\sb T$ let 
$top(A)/sim^1_A=\{T_i:i\in I_A\}$, by our  
assumption, for every $l<\om$ there is an initial segment $A\sb T$ \st 
$|I_A|>l$. Choose a large enough $l$ (see below) and a corresponding $A$ 
and for every $i\in I_A$ pick $x_i\in T_i$.

\\If\/ $l$ is larger than the number of possible theories 
($=|T_{n,l(\bbp)}|$) then there are 
$i\not= j\in I_A$ \st $Th^n(T_i;x_i,\bb P)=Th^n(T_j;x_j,\bb P)$ and let's 
assume that we have chosen such an $l$. Now let 
$\bb R_1=\{x_i\}\cup\{x_i,x_j\}\cup\bb P$ and 
$\bb R_2=\{x_j\}\cup\{x_i,x_j\}\cup\bb P$.
Apply 1.12: clearly 
$$Th^m(T_{\le A};\bb R_1)=Th^m(T_{\le A};\bb R_2)=
Th^m(T_{\le A};\em,\em,\bb P)$$ 
\\and easily 
$$Th^m(I_A;\bb Q^A(n,\bb R_1)=Th^m(I_A;\bb Q^A(n,\bb R_2)$$ 
\\but by 1.12 these theories determine $Th^n(T;x_i,\{x_i,x_j\},\bb P)$ and 
$Th^n(T;x_j,\{x_i,x_j\},\bb P)$ hence 
$$T\models \varphi(x_i,\{x_i,x_j\},\bb P) \iff 
T\models \varphi(x_j,\{x_i,x_j\},\bb P)$$ 
\\a contradiction. \vv \qed
\pt Claim 4.4. If\/ $T$ is a wild tree and (1)$(ii)$ of 4.2 holds 
then no monadic 
formula \ $\varphi(x,X,\bb P)$ defines a choice function on $T$. \vv
\pd Proof. Let $B\sb T$ be a branch that embeds $\Q$. We will apply 1.13 and 
``translate'' the choice function on $T$ to a choice function on $B$ but by 
2.4 there is no definable choice function on $B$.  

\\So assume that $\varphi(x,X,\bbp)$ defines a choice function on $T$ and is 
determined by $Th^n(T;x,X,\bbp)$. By 1.13 there is an $m<\om$, a chain $B'$ 
with $(B,\seg)\sb(B',\seg)$ and a sequence of parameters $\bb Q\sb B'$ \st 
from $Th^m(B';\bb Q)$ we can compute $Th^n(T;\bbp)$. Define, for 
$\et\seg\nu\in B$, $f(\et,\nu)=Th^m(B';\bb Q)\red_{[\et,\nu)}$.  
$f$ is an additive colouring hence by 2.3 there is $X=\{\et_i\}_{i\in\Z}$, of 
order type $\Z$, homogeneous with respect to $f$. As in the proof of 2.4 we 
have: 
$$i,j\in\Z \imp Th^m(B';\et_i,X,\bb Q)=Th^m(B';\et_j,X,\bb Q)$$
and (by the `moreover' clause in 1.13) this implies 
$$i,j\in\Z \imp Th^n(T;\et_i,X,\bbp)=Th^n(T;\et_j,X,\bbp).$$
Hence
$$i,j\in\Z \imp 
[T\models\varphi(\et_i,X,\bbp)\iff T\models\varphi(\et_j,X,\bbp)]$$
and this contradicts ``$\varphi$ chooses an element from $X$''.
\vv \qed
\pt Claim 4.5. If\/ $T$ is a wild tree and (1)$(iii)$ of 4.1 holds then no 
monadic formula \ $\varphi(x,X,\bb P)$ defines a choice function on $T$. \vv
\pd Proof. Similar to the previous proof. 

\\By (1)$(iii)$ for every $m<\om$ there is a branch $B\sb T$ with 
$\hd(B)>m$. Use 1.13, 3.10 and the proof of 3.9 to find, for a suitable 
branch $B$, a homogeneous subset that contradicts the assumption that
$\varphi(x,X,\bb P)$ defines a choice function on $T$. 

\\The details are left to the reader.
\vv \qed
\\We conclude
\pt Theorem 4.6.  $T$ is a wild tree \ $\imp$ \  $T$ does not have a 
monadically definable choice function. Moreover, every candidate fails to  
choose from either linearily ordered subsets (4.4, 4.5) or anti-chains (4.3).
\vv \qed
\sec
%
%
%
%
%
\hh 5. Tame trees \par
By \gu\ $\bina$ does not have a definable choice function. To know if a tame 
tree $T$ has a definable choice function we just have to ask if there is an 
embedding of\/ $f\colon\bina\too T$. If such an embedding exists we use \gu\ 
to show that $T$ does not have one, if not, $T$ has even a definable well 
ordering.
\pt Claim 5.1. Let $T$ be a tree and $F\colon\bina\too T$ be a tree embedding.
Then no monadic formula \ $\varphi(x,X,\bb P)$ defines a choice function on 
$T$. \vv
\pd Proof. We will use \gu\ 1.15 and the notations of 1.14. First, we may 
assume w.l.o.g that $T$ has a root (adding a root will not effect the 
existence of a choice function) and that $F(root(\bina))=root(T)$. 
Now apply the proof in $\S$5 of \gu. From the proof there we learn that 
for every $\bb Q\sb\bina$ and $m<\om$ there is an infinite anti-chain 
$Y\sb\bina$ \st for every $y\in Y$ there is $y^*\not=y\in Y$ with 
$Th^m(Bush_{\bina}(Y);y,\bb Q)=Th^m(Bush_{\bina}(Y);y^*,Y,\bb Q)$. 
In our context ($F''(\bina)=S\sb T$) the result has the form:

$(*)$ for every $\bb Q\sb S$ and $m<\om$ there is an infinite anti-chain
$Y\sb S$ \st for every $y\in Y$ there is $y^*\not=y\in Y$ with 
$Th^m(Bush_S(Y);y,\bb Q)=Th^m(Bush_S(Y);y^*,\bb Q)$.

\\Let $\varphi(x,X,\bb P)$ be a candidate for a definition of a choice 
function on $T$ and suppose $Th^n(T;x,X,\bbp)$ decides $\varphi$. 
Let $m<\om$ and $\bb Q=\bb Q(n,\bbp)$ be as in 1.15 and $Y\sb S$ be the 
anti-chain from $(*)$. 
Suppose $T\models\varphi(y,Y,\bb P)$, by $(*)$ we have $y^*\in Y$ as in there. 
Now $Th^m(Bush_S(Y);y,\bb Q)=Th^m(Bush_S(Y);y^*,\bb Q)$ and by 1.15 
$$Th^n(T;y,Y,\bbp)=Th^n(T;y^*,Y,\bbp)$$
hence
$$T\models\varphi(y,Y,\bb P)\iff T\models\varphi(y^*,Y,\bb P)$$
hence $\varphi$ fails to define a choice function on $T$.
\vv \qed
\pd \dd 5.2. Let $T$ be a tree. For $\et\in T$ we define by recursion a rank 
function $rk(\et)$ by:  

\\$rk(\et)\ge{\al+1} \iff$ there are $\nu_1,\nu_2\in T$ with $\et\seg\nu_1$ and 
$\et\seg\nu_2$ \st $\nu_1,\nu_2$ are incomparable in $T$ and 
$rk(\nu_1),rk(\nu_2)\ge\al$

\\If\/ $rk(\et)$ is not defined we stipulate $rk(\et)=\infty$. \vv
\pt Fact 5.3.  (1) $\et\seg\nu\in T \ \imp \ rk(\nu)\le rk(\et)$ where 
$\le$ has the obvious meaning.

\\(2) $\bina$ is not embeddable in a tree $T$ \ $\iff$ \ for every $\et\in T$, 
$rk(\et)\not=\infty$ \vv
\pt Lemma 5.4.  Let $T$ be a tame tree. If\/ $\bina$ is not embeddable in 
$T$ then there are $\bb Q\sb T$ and a monadic formula $\varphi(x,y,\bb Q)$ 
that defines a well ordering of\/ $T$. \vv
\pd Proof. Assume $T$ is $(n^*,k^*)$ tame, recall definitions 4.1 and 4.2 and 
remember that for every $x \in T$, $rk(x)$ is well defined (i.e. $<\infty$).
We will partition $T$ into a disjoint union of sub-branches, 
indexed by the nodes of a well founded tree $\Ga$ and reduce the problem of a 
well ordering of\/ $T$ to a problem of a well ordering of\/ $\Ga$.

Step 1. Define by induction on $\al$ a set $\Ga_\al\sb\thinspace^\al Ord$ 
(this is a our set of indices), for every $\et\in\Ga_\al$ define a tree 
$T_\et\sb T$ and a branch $A_\et\sb T_\et$.

\\$\al=0$ : $\Ga_0$ is $\{\lan \ran\}$, $T_{\lan \ran}$ is $T$ and 
$A_{\lan \ran}$ is a branch (i.e. a maximal linearily ordered subset) of\/ $T$.

\\$\al=1$ : Look at $(T\setminus A_{\lan \ran})/\sim^1_{A_{\lan \ran}}$, it's 
a disjoint union of trees and name it $\lan T_{\lan i\ran}:i<i^* \ran$,
let $\Ga_1:=\{\lan i\ran:i<i^*\}$ and for every $\lan i\ran\in\Ga_1$ let 
$A_{\lan i\ran}$ be a branch of\/ $T_{\lan i\ran}$.

\\$\al=\be+1$ : For $\et\in\Ga_\be$ denote 
$(T_\et\setminus A_\et)/\sim^1_{A_\et}$ by $\{T_{\ee\et,i}:i<i_\et\}$, let  
$\Ga_\al=\{\ee\et,i:\et\in\Ga_\be,\ i<i_\et\}$ and choose $A_{\ee\et,i}$ to 
be a branch of\/ $A_{\ee\et,i}$.

\\$\al$ limit: Let $\Ga_\al=
\{\et\in\thinspace^\al Ord : \wedge_{\be<\al}\et\red_\be\in\Ga_\be,\ 
\wedge_{\be<\al}T_{\et\red_\be}\not=\emptyset\}$, 
let for $\et\in\Ga_\al$ $T_\et=\cap_{\be<\al}T_{\et\red_\be}$ and 
$A_\et$ a branch of\/ $T_\et$. ($T_\et$ may be empty).

\\Now, at some stage $\al\le|T|^+$ we have $\Ga_\al=\emptyset$ and let 
$\Ga=\cup_{\be<\al}\Ga_\be$. Clearly $\{A_\et:\et\in\Ga\}$ is a partition of 
$T$ into disjoint sub-branches.

\\Notation: having two trees $T$ and $\Ga$, to avoid confusion, we use 
$x,y,s,t$ for nodes of\/ $T$ and $\et,\nu,\si$ for nodes of\/ $\Ga$.

Step 2. We want to show that $\Ga_\om=\emptyset$ hence $\Ga$ is a well 
founded tree. Note that we made no restrictions on the choice of 
the $A_\et$'s and we add one now in order to make the above statement true.
Let $\ee\et,i\in\Ga$ define $A_{\et,i}$ to be the sub-branch 
$\{t\in A_\et:(\all s\in A_{\ee\et,i})[rk(t)\le rk(s)]\}$ and 
$\ga_{\et,i}$ to be $rk(t)$ for some $t\in A_{\et,i}$. By 5.5(1) and the 
inexistence of a stricly decreasing sequence of ordinals, 
$A_{\et,i}\not=\emptyset$ and $\ga_{\et,i}$ is well defined. Note also that 
$s\in A_{\ee\et,i} \imp rk(s)\le\ga_{\et,i}$. 

\\\underbar {Proviso}: For every $\et\in\Ga$ and $i<i_\et$ the sub-branch 
$A_{\ee\et,i}$ contains every $s\in T_{\ee\et,i}$ with $rk(s)=\ga_{\et,i}$. 

\\Following this we claim: ``$\Ga$ does not contain an infinite, stricly 
increasing sequence''.  Otherwise let $\{\et_i\}_{i<\om}$ be one, and 
choose $s_n\in A_{\et_n,\et_{n+1}(n)}$ (so $s_n\in A_{\et_n}$). Clearly 
$rk(s_n)\ge rk(s_{n+1})$ and by the proviso we get 
$$rk(s_n)=rk(s_{n+1})\imp rk(s_{n+1})>rk(s_{n+2})$$
\\therefore $\{rk(s_n)\}_{n<\om}$ contains an infinite, stricly decreasing 
sequence of ordinals which is absurd.

Step 3. Next we want to make ``$x$ and $y$ belong to the same $A_\et$'' 
definable. 

\\For each $\et\in\Ga$ choose $s_\et\in A_\et$, and let $Q\sb T$ be the set
of representatives.
Let $h\colon T\to\{d_0,\ldots,d_{n^*-1}\}$ be a colouring that satisfies:
$h\red_{A_{\lan \ran}}=d_0$ and for every $\ee\et,i\in\Ga$, 
$h\red_{A_{\ee\et,i}}$ is constant and, when $j<i$ and 
$s_{\ee\et,j}\sim^0_{A_\et}s_{\ee\et,j}$ we have  
$h\red_{A_{\ee\et,i}}\not=h\red_{A_{\ee\et,j}}$. This can be done as $T$ is 
$(n^*,d^*)$ tame.

\\Using the parameters $D_0,\ldots,D_{n^*-1}$ ($x\in D_i$ iff\/ $h(x)=d_i$),
we can define $\vee_\et x,y\in A_\et$ by ``$x,y$ are comparable and the 
sub-branch $[x,y]$ (or $[y,x]$) has a constant colour''.

Step 4. As every $A_\et$ has Hausdorff degree at most $k^*$, we can define a 
well ordering of it using parameters $P^\et_1,\ldots,P^\et_{k^*}$ and by 
taking $\bb P$ to be the (disjoint) union of the $\bb P^\et$'s we can define 
a partial ordering on $T$ which well orders every $A_\et$.

\\By our construction $\et\seg\nu$ if and only if there is an element in 
$A_\nu$ that `breaks' $A_\et$ i.e. is above a proper initial segment of 
$A_\et$. (Caution, if\/ $T$ does not have a root this may not be the case for  
$\lan \ran$ and a $<n^*$ number of\/ $\lan i\ran$'s and we may need parameters 
for expressing that). Therefore, as by step 3 ``being in the same $A_\et$'' 
is definable, we can define a partial order on the sub-branches $A_\et$ 
(or the representatives $s_\et$) by $\et\seg\nu\imp A_\et\le A_\nu$.  

\\Next, note that ``$\nu$ is an immediate successor of\/ $\et$ in $\Ga$'' is 
definable as a relation between $s_\nu$ and $s_\et$ hence the set 
$A_\et^+ := A_\et\cup\{s_{\ee\et,i}\}$ is definable from $s_\et$.  Now 
the order on $A_\et$ induces an order on $\{s_{\ee\et,i}/\sim^0_{A_\et}\}$ 
which is can be embedded in the complition of\/ $A_\et$ hence has Hdeg$\le k^*$. 
Using additional parameters $Q^\et_1,\ldots,Q^\et_{k^*}$, we have a definable 
well ordering on $\{s_{\ee\et,i}/\sim^0_{A_\et}\}$. As for the ordering on 
each $\sim^1_{A_\et}$ equivalence class (finite with $\le n^*$ elements),
define it by their colours (i.e. the element with the smaller colour is the 
smaller according to the order).

\\Using $\bb D$, $\bb P$, $Q$ and $\bb Q=\cup_\et\bb Q^\et$ we can define a 
partial ordering which well orders each $A_\et^+$ in such a way that every
$x\in A_\et$ is smaller then every $s_{\ee\et,i}$.

\\Summing up we can define (using the above parameters) a partial order 
on subsets of\/ $T$ that well orders each $A_\et$, orders sub-branches 
$A_\et$, $A_\nu$ when the indices are comparable in $\Ga$ and well 
orders all the ``immediate successors'' sub-branches of a sub-branch $A_\et$.

Step 5.  The well ordering of\/ $T$ will be defined by $x<y \iff$

\\a) $x$ and $y$ belong to the same $A_\et$ and $x<y$ by the well order on 
$A_\et$;  \ or

\\b) $x\in A_\et$, $y\in A_\nu$ and $\et\seg\nu$;  \ or

\\c) $x\in A_\et$, $y\in A_\nu$,  $\si=\et\wedge\nu$ in $\Ga$ (defined as a 
relation between sub-branches), $\ee\si,i\seg\et$, $\ee\si,j\seg\nu$  
and $s_{\ee\si,i}<s_{\ee\si,j}$ in the order of\/ $A^+_\si$.

\\Note, that $<$ is a linear order on $T$ and every $A_\et$ is a convex 
and well ordered sub-chain. Moreover $<$ is a linear order on $\Ga$ and the 
order on the $s_\et$'s is isomorphic to a lexicographic order on $\Ga$.

\\Why is the above (which is clearly definable with our parameters) a well 
order?  Because of the above note and because a lexicographic ordering of a  
well founded tree is a well order, provided that immediate successors are 
well ordered. In detail, assume $X=\{x_i\}_{i<\om}$ is a stricly decreasing 
sequence of elements of\/ $T$.  Let $\et_i$ be the unique node in $\Ga$ \st 
$x_i\in A_{\et_i}$ and by the above note w.l.o.g 
$i\not=j\imp\et_i\not=\et_j$.  By the well foundedness of\/ $\Ga$ and clause 
(b) we may also assume w.l.o.g that the $\et_i$'s form an anti-chain in 
$\Ga$. Look at  $\nu_i:=\et_1\wedge\et_i$  which is constant for infinitely 
many $i$'s and w.l.o.g equals to $\nu$ for every $i$.  Ask:

\\$(*)$ is there is an infinite $B\sb\om$ \st 
$i,j\in B\imp x_i\sim^0_{A_\nu}x_j$\ ?

\\If this occurs we have $\nu_1\not=\nu$ with $\nu\seg\nu_1$ \st for 
some infinite $B'\sb B\sb\om$ we have $i\in B'\imp \nu_1\seg\et_i$. 
(use the fact that $\sim^1_{A_\nu}$ is finite).
W.l.o.g $B'=\om$ and we may ask if $(*)$ holds for $\nu_1$. 
Eventually, since $\Ga$ does not have an infinite branch, we will have a 
negative answer to $(*)$. We can conclude that w.l.o.g there is 
$\nu\in\Ga$ \st $i\not=j\imp x_i\not\sim^0_{A_\nu}x_j$ i.e. the $x_i$'s 
``break'' $A_\nu$ in ``different places''.  

\\Define now $\nu_i$ to be the unique immediate successor of\/ $\nu$ \st 
$\nu_i\seg\et_i$. The set $S=\{s_{\nu_i}\}_{i<\om }\sb A^+_\nu$ is well 
ordered by the well ordering on $A^+_\nu$ and by clause (c) in the definition 
of\/ $<$,  $x_i>x_j\iff\nu_i>\nu_j$ so $S$ is an infinite stricly decreasing 
subset of\/ $A^+_\nu$ -- a contradiction.

\\This finishes the proof that there is a definable well order of\/ $T$.
\vv \qed

\\Finally we can conclude: 
\pt Theorem 5.5. Let $T$ be a tree. If\/ $T$ is wild or $T$ embeds $\bina$ then 
there is no definable choice function on $T$ (by a monadic formula with 
parameters). If\/ $T$ is tame and does not embed $\bina$ then there even 
a definable well ordering of the elements of\/ $T$ by a monadic formula (with 
parameters) $\varphi(x,y,\bb P)$.
\vv \qed

\\As mentioned in the introduction, a tree is tame [wild] [embeds $\bina$]
if and only if it's completion is tame [wild] [embeds $\bina$]. 
Moreover looking at the proofs of 4.3, 4.4, 4.5 and 5.1 we note that the 
counter-examples for the choice function problem are either anti-chains or 
linearily ordered subsets of\/ $T$. We conclude:

\pt Conclusion 5.6. Let $T$ be a tree and $T'$ be it's completion. Then the 
following are equivalent:

\\a) For some $n,l<\om$, for every anti-chain/branch $A$ of\/ $T$ there is a 
monadic formula $\varphi_A(x,X,\bb P_A)$ with quantifier depth $\le n$ and 
$\le l$ parameters from $T$, that defines a choice function from non empty 
subsets of\/ $A$. 

\\b) There is a monadic formula, with parameters, $\psi(x,y,\bb P)$ that 
defines a well ordering of the elements of\/ $T$.

\\c) There is a monadic formula, with parameters, $\psi'(x,y,\bb P')$ that 
defines a well ordering of the elements of\/ $T'$.
\vv \qed
\sec
\eject
%
%

\font\ba=cmr8
\font\bs=cmbxti10
\font\bib=cmtt12
\centerline{\bib REFERENCES}  \medbreak
\\ \ {\bf [BL]} \ J.R. B{\ba \"UCHI} and L.H. L{\ba ANDWEBER}, \ 
{\sl Solving sequential conditions by finite-state strategies}, \ 
{\bs Transactions of the American Mathematical Society}, 
\ vol. 138 (1969), pp. 295--311.   \vv
\\ \ {\bf [Ha]} \ F. H{\ba AUSDORFF}, \ 
{\sl Grundz\"uge einer Theorie der geordnetn Mengen}, \ 
{\bs Math. Ann.}, 
\ vol. 65 (1908), pp. 435--505.   \vv
\\ \ {\bf [Gu]} \ Y. G{\ba UREVICH},  \ 
{\sl Monadic Second--order Theories}, \ 
{\bs Model Theoretic Logics}, \ 
(J. Barwise and S. Feferman, editors), 
\ Springer--Verlag, Berlin 1985, pp. 479--506.  \vv
\\ \ {\bf [GuSh]} \ Y. G{\ba UREVICH} and S. S{\ba HELAH}, \ 
{\sl Rabin's Uniformization Problem}, \ 
{\bs The Journal of Symbolic Logic}, 
\ vol. 48 (1983), pp. 1105--1119.   \vv
\\ \ {\bf [Ra]} \ M.O. R{\ba ABIN}, \ 
{\sl Decidability of second-order theories and automata on infinite trees}, \ 
{\bs Transactions of the American Mathematical Society}, 
\ vol. 141 (1969), pp. 1--35.   \vv
\\ \ {\bf [Sh]} \ S. S{\ba HELAH}, \ 
{\sl The monadic Theory of Order}, \ 
{\bs Annals of Mathematics}, 
\ ser. 2, vol. 102 (1975), pp. 379--419.   \vv
\end